\documentclass[twoside,11pt]{amsart}

\usepackage{amsmath}
\usepackage{amsfonts}
\usepackage{amssymb}
\usepackage{latexsym}
\usepackage{epsfig}

\usepackage{float}
\restylefloat{figure}
\usepackage{subfigure}

\usepackage{epic}

\usepackage{graphicx}
\usepackage{graphics}
\usepackage{moreverb}
\usepackage{multicol}

\newtheorem{theorem}{Theorem}[section]
\newtheorem{lemma}{Lemma}[section]

\newtheorem{nota}{Remark}[section]
\newtheorem{cor}{Corollary}[section]
\newcommand{\disp}{\displaystyle}

\newcommand{\barint}{\hbox{$\int$\kern-0.75\intwidth
\vrule width 0.5\intwidth height 2.4pt depth
-2pt\kern0.25\intwidth}}
\newlength\intwidth
\setbox0=\hbox{$\int$} \intwidth=\wd0

\newcommand\avint{\hbox{\hbox{$\displaystyle \int$}\hbox{\kern-.9em{$-$}}}}

\newcommand\smavint{\hbox{\hbox{$\int$}\hbox{\kern-.75em{$-$}}}}

\def\XXint#1#2#3{{\setbox0=\hbox{$#1{#2#3}{\int}$}
\vcenter{\hbox{$#2#3$}}\kern-.51\wd0}}

\newcommand{\R}{\mbox{\rm\bf R}}
\newcommand{\eps}{\epsilon}
\newcommand{\be}{\begin{eqnarray}}
\newcommand{\ee}{\end{eqnarray}}

\newcommand{\ben}{\begin{eqnarray*}}
\newcommand{\een}{\end{eqnarray*}}

\newcommand{\srho}{\sqrt{\rho}}

\newcommand{\dac}{\backprime}

\numberwithin{equation}{section}

\begin{document}

\title[On the  global regularity of sub-critical  Euler-Poisson equations]
{On the  global regularity of\\sub-critical  Euler-Poisson equations with pressure}

\author[Eitan Tadmor]{Eitan Tadmor}
\address[Eitan Tadmor]{\newline
        Department of Mathematics,  Institute for Physical Science and Technology\newline
        and Center of Scientific Computation And Mathematical Modeling (CSCAMM)\newline
        University of Maryland\newline
        College Park, MD 20742 USA}
 \email[]{tadmor@cscamm.umd.edu}
\urladdr{http://www.cscamm.umd.edu/\~{}tadmor}
\author[Dongming Wei]{Dongming Wei}
\address[Dongming Wei]{\newline
    Department of Mathematics\newline
    Center of Scientific Computation And Mathematical Modeling (CSCAMM)\newline
        University of Maryland\newline
        College Park, MD 20742 USA}
\email []{dwei@cscamm.umd.edu}
\urladdr{http://www.cscamm.umd.edu/people/students/dwei.htm}

\date\today

\subjclass{35Q35, 35B30.}
\keywords{Euler-Poisson equations, Riemann Invariants, Critical thresholds, Global regularity.}

\thanks{\textbf{Acknowledgment.} Research was supported in part by NSF grant 04-07704 and ONR grant N00014-91-J-1076.}


\begin{abstract}

We prove that the one-dimensional Euler-Poisson system driven by the
Poisson forcing together with the usual $\gamma$-law
pressure, $\gamma \geq 1$, admits  global solutions for a large class
of initial data. Thus, the Poisson forcing regularizes the generic
finite-time breakdown in the $2 \times 2$ p-system. Global
regularity is shown to depend on whether or not the  initial configuration
of the Riemann invariants and density crosses an intrinsic critical
threshold.

\end{abstract}
\maketitle

\section{Introduction}\label{Intro}
\setcounter{figure}{0}
\setcounter{equation}{0}

It is well known that the systems of Euler
equations for compressible flows can and will breakdown at a finite
time even if the initial data are smooth. A prototype example for
such systems is provided by the $2 \times 2$ system of isentropic
gas dynamics
\begin{equation}\label{ES} \left\{ \begin{array}{l}
 \rho_t+(\rho u)_x=0\, \\
(\rho u)_t+(\rho u^2)_x=-p_x \,, \\
\end{array} \right.
\end{equation}
where the pressure $p=p(\rho)$ is given by the usual $\gamma$-law,
$p(\rho)=A\rho^{\gamma}$. By using the method introduced in
\cite{La64} to deal with pairs of conservation laws, it can be
shown that (\ref{ES}) will lose the $C^1$-smoothness due to the
appearance of shock discontinuities unless its two Riemann
invariants are nondecreasing. Thus, the finite time breakdown of
(\ref{ES}) is generic in the sense that it holds for all but a
``small set" of initial data.

On the other hand, if we replace the pressure by Poisson forcing,
then we arrive at the system of Euler-Poisson equations
\begin{equation}\label{PEP} \left\{ \begin{array}{l}
 \rho_t+(\rho u)_x=0\,, \\
(\rho u)_t+(\rho u^2)_x =-k\rho\varphi_x \, \qquad k>0 \,, \\
\end{array} \right.
\end{equation}
subject to initial data $(u_0\,, \rho_0>0)$\,. Here $\varphi=\varphi(\rho)$
is the potential, which is dictated by the (one-dimensional) Poisson
equation,\,$\varphi_{xx}=-\rho$\,. In this case, there is a
``large set" of initial configurations which yield global smooth
solutions. More precisely, \cite{ELT01} have shown that (\ref{PEP})
admits a global smooth solution if and only if
\begin{equation}\label{epct}
u_{0x}(x)>-\sqrt{2k\rho_0(x)}\,.
\end{equation}
Thus, following the terminology of \cite{LT02}, the curve $u_{0x}+\sqrt{2k\rho_0}=0$ is a  ``critical threshold" in configuration space which separates between initial
configurations leading to finite time breakdown and a ``large
set" of sub-critical initial configurations which yield global smooth solutions. In particular, (\ref{epct}) allows negative velocity gradients (-- depending on the local amplitude of the density), which
otherwise are excluded in the case of inviscid Burgers equations, corresponding to $k=0$.

In this paper we turn our attention to the full Euler-Poisson
equations driven by both -- pressure and Poisson forcing,
\begin{equation}\label{FEP} \left\{ \begin{array}{l}
\rho_t+(\rho u)_x=0\,, \\
(\rho u)_t+(\rho u^2)_x =-p(\rho)_x-k\rho\varphi_x, \quad k>0, \\
-\varphi_{xx}=\rho\,.
\end{array} \right.
\end{equation}

These equations govern different phenomena, ranging from
the largest scale of e.g., the evolution  gravitational collapse  in stars, to applications in the smallest scale of e.g., semi-conductors.
There is a considerable amount of literature available on the local and global
behavior of Euler-Poisson and related problems. Consult \cite{Ma86} for
local existence in the small $H^s$-neighborhood of a steady state of
self-gravitating stars,
\cite{CW96} for global existence of weak solutions with geometrical
symmetry, \cite{Gu98} for global existence for 3-D irrotational flow,
\cite{MN95} for isentropic case, and \cite{JR00} \cite{PRV95} for
isothermal case. Consult \cite{Pe90} \cite{MP90} \cite{Si85} \cite{En96}
\cite{WC98}, \cite{BW98} and in particular, \cite{En96} (more about that below),
for non-existence results and singularity formation.
The question of global smoothness vs. finite breakdown was studied in a recent
series of works of Engelberg, Liu and Tadmor, in terms of a critical threshold phenomena for
1-D ``pressure-less" Euler-Poisson equations, \cite{ELT01} and
 2-D restricted Euler-Poisson equations, \cite{LT02, LT03}.

The natural question that arises in the present context of full Euler-Poisson equations (\ref{FEP}) is whether the pressure  enforces a generic finite time
breakdown or, whether the presence of Poisson forcing preserves global regularity for a ``large" set of initial configurations.
We answer this question of ``competition" between pressure and Poisson forcing, proving that the Euler-Poisson equations (\ref{FEP}) with $\gamma \geq 1$ admit global smooth solutions for a ``large set" of sub-critical initial data such that
\begin{equation}\label{eq:gonea}
u_{0x}(x) > -K_0\sqrt{\rho_0(x)} +
\sqrt{A\gamma}\frac{|\rho_{0x}(x)|}{\rho_0(x)^{\frac{3-\gamma}{2}}},
\qquad \gamma \geq 1.
\end{equation}
Here, $K_0$ is a constant depending on $k, \gamma$ and the initial data. In the particular (and important) case of isothermal equations, $\gamma=1$, we have $K_0=\sqrt{2k}$ and (\ref{eq:gonea}) amounts to a sharp critical threshold,
\begin{equation}\label{eq:goneb}
u_{0x}(x) \geq -\sqrt{2k\rho_0(x)} +
\sqrt{A}\frac{|\rho_{0x}(x)|}{\rho_0(x)}, \qquad \gamma=1.
\end{equation}
The inequalities (\ref{eq:gonea}),(\ref{eq:goneb}) quantify the competition between the destabilizing pressure effects, as the range of sub-critical initial configurations shrinks with the growth of the amplitude of the pressure, $A$, while the stabilizing effect of the Poisson forcing increases the sub-critical range with a growing $k$.
In particular, (\ref{eq:goneb}) with $A=0$  recovers the pressure-free critical threshold (\ref{epct}).

Formation of singularities and global regularity of (\ref{FEP}) were addressed earlier by Engelberg in \cite{En96}. His results show finite-time break-down if  
$u_{0x}(x) - \sqrt{A\gamma}|\rho_{0x}(x)|\rho_0(x)^{\frac{\gamma-3}{2}}$ is ``...sufficiently negative at some point". Our contribution here is to \emph{quantify} the critical threshold behind this asymptotic statement. To fully appreciate this quantified threshold, we turn to the converse statement in \cite[Theorem 2]{En96}: it asserts the global regularity of (\ref{FEP}) for a class of initial data such that $u_{0x}(x) - \sqrt{A\gamma}|\rho_{0x}(x)|\rho_0(x)^{\frac{\gamma-3}{2}}>0$.
It is a ``non-generic" class (in the sense of requiring both Riemann invariants at $t=0$ to be \emph{globally} increasing). 
In fact, by (\ref{eq:gonea}) one has a \emph{negative} threshold, $-K_0\sqrt{\rho_0}$, implying the existence of a ``large" class of sub-critical initial data with global regularity.
 
The paper is organized as follows. In section \ref{RV}, we reformulate the
system (\ref{FEP}) with its Riemann invariants as a preparation for
the analysis carried out in sections \ref{CT} and \ref{FTB}. In section \ref{CT}, we prove our main results, providing sufficient conditions for ``large sets" of sub-critical initial configurations which yield global smooth solution.
In section \ref{FTB}, we  give examples of  finite time
breakdown for super-critical initial data. Combining our results in sections \ref{CT} and \ref{FTB}, they confirm the existence of a critical threshold phenomena for the full Euler-Poisson equations (\ref{FEP}).

\section{Riemann Invariants} \label{RV}

\setcounter{figure}{0}
\setcounter{equation}{0}

\subsection{The Euler-Poisson equations with $\gamma$-law pressure: $\gamma >1$}
We begin by rewriting the Euler-Poisson equations (\ref{FEP}) as a
first order quasilinear system
\begin{equation} \left( \begin{array}{l} \rho \\ u \end{array}
\right)_t + J \left( \begin{array}{l} \rho \\ u \end{array}
\right) _x=\left( \begin{array}{l} 0 \\ -k\varphi_x
\end{array} \right) \,,
\end{equation}
where the Jacobian $J:= \left( \begin{array}{ll} u & \rho \\
A\gamma\rho^{\gamma-2} & u \end{array} \right)$ has two different
eigenvalues
\[
\disp{\lambda:=u-\sqrt{A\gamma}\rho^{\frac{\gamma-1}{2}}
<\mu:=u+\sqrt{A\gamma}\rho^{\frac{\gamma-1}{2}}\,.}
\]
and let $R$ and $S$ denote the Riemann invariants of the
corresponding Euler system (\ref{ES})
\begin{equation}\label{eq:RS}
R:=u-\frac{2\sqrt{A\gamma}}{\gamma-1}\rho^{\frac{\gamma-1}{2}}
 \quad \mbox{ and }
\quad S:=u+\frac{2\sqrt{A\gamma}}{\gamma-1}\rho^{\frac{\gamma-1}{2}}\,.
\end{equation}
They satisfy the coupled system of equations,
\begin{subequations}\label{eq:RI}
\begin{align}\label{RIa}
R_t+\lambda R_{x}=-k\varphi_x\,, \\
\label{RIb} S_t+\mu S_{x}=-k\varphi_x\,,
\end{align}
\end{subequations}
coupled through the Poisson equation $-\phi_{xx}=\rho\,.$
If we set $r:=R_{x}\,,$ $s:=S_{x}\,$ then upon differentiation of (\ref{eq:RI})
we get
\begin{subequations}\label{eq:pair}
\begin{align}
r_t+\lambda r_{x}+\lambda_{S}rs+\lambda_{R}r^2=k\rho\,, \\
s_t+\mu s_{x}+\mu_{S}s^2+\mu_{R}rs=k\rho\,.
\end{align}
\end{subequations}

\noindent
Next, we observe that $\lambda=\frac{R+S}{2}-\frac{\gamma-1}{4}(S-R)$ and
$\mu=\frac{R+S}{2}+\frac{\gamma-1}{4}(S-R)$. Hence, expressed in terms of
${\displaystyle \theta:=\frac{\gamma-1}{2}}$, we have for $\gamma \geq 1$, 
\[
\lambda_R=\mu_S=\frac{1+\theta}{2}\quad \mbox{and}\quad \lambda_S=\mu_R=\frac{1-\theta}{2},
\qquad \theta:=\frac{\gamma-1}{2} \geq 0,
\]
 and the  pair of equations (\ref{eq:pair}) is recast into the form
\begin{subequations}
\begin{align} \label{r2}
r^{\dac}+\frac{1+\theta}{2} r^2+\frac{1-\theta}{2} rs&=k\rho,\\
\label{s2} {s}^{\prime} +\frac{1+\theta}{2} s^2+\frac{1-\theta}{2} rs&=k\rho\,.
\end{align}
\end{subequations}
Here and below $\{\}^{\dac}:=\partial_t+\lambda\partial_x$
and $\{\}':=\partial_t+\mu\partial_x$ denote differentiation along the $\lambda$
and $\mu$ particle paths,
\[
\Gamma_{\lambda}:=\big\{(x,t)\,\big|\,
\disp \dot{x}(t)=\lambda\big(\rho(x,t),u(x,t)\big)\big\},
\quad
\Gamma_{\mu}:=\big\{(x,t)\,\big|\,
\disp\dot{x}(t)=\mu\big(\rho(x,t),u(x,t)\big) \big\}.
\]
To continue, we rewrite the equation for $\rho$ as
\begin{equation}\label{eq:rho}
(\rho_t+\lambda
\rho_x)+\disp{\frac{\mu-\lambda}{2}}\rho_x+\rho{\disp\frac{s+r}{2}}=0\,.
\end{equation}
Since
$\disp s-r=S_x-R_x=2\sqrt{A\gamma}\rho^{\frac{\gamma-3}{2}}\rho_x$\,,
it enables us to express
$\disp \frac{\mu-\lambda}{2}\rho_x=\sqrt{A\gamma}\rho^{\frac{\gamma-1}{2}}\rho_x=\rho\frac{s-r}{2}\,, $ so that the $\rho$ equation (\ref{eq:rho}) can be
written along the $\lambda$ particle path as
$\rho ^{\dac}+\rho s=0$.
Similarly, it can be written along the $\mu$ particle path
as ${\rho}^{\prime}+\rho r=0$.
Assembling the above equations together,
we arrive at the  following system governing $r$, $s$ and $\rho$,
\begin{subequations}\label{pair:sys}
\begin{equation}\label{sys}
\left\{ \begin{array}{rl}
{\displaystyle r^{\dac} +\frac{1+\theta}{2} r^2+\frac{1-\theta}{2} rs} & =k\rho \,, \\
\rho^{ \dac} +\rho s& =0 \,,
\end{array} \right.
\end{equation}
and
\begin{equation}\label{sysa}
\left\{  \begin{array}{rl}
{\displaystyle {s}^{\prime}+\frac{1+\theta}{2}
s^2+\frac{1-\theta}{2} rs} & = k\rho \,, \\
{\rho}^{\prime}+\rho r& =0.
\end{array} \right.
\end{equation}
\end{subequations}

\noindent
Finally, we use the integration factors $1/\srho$ and $r/2\rho\srho$ in the first and second equations of each pair in (\ref{pair:sys}), to conclude
\begin{subequations}\label{pair:good}
\begin{align}\label{eq:gooda}
{\displaystyle \bigg(\frac{r}{\srho}\bigg)^{\dac} +\frac{1+\theta}{2}
\frac{r^2}{\srho}
- \frac{\theta}{2}\frac{rs}{\srho}}&=k\srho, \\
\label{eq:goodb}
{\displaystyle \bigg(\frac{s}{\srho}\bigg)^{\prime}+\frac{1+\theta}{2}
\frac{s^2}{\srho}- \frac{\theta}{2}\frac{rs}{\srho}} &=k\srho.
\end{align}
\end{subequations}

\subsection{The isothermal case $\gamma=1$}\label{subsec:iso}
In this case, the two eigenvalues are $\lambda=u-
\sqrt{A}<\mu=u+\sqrt{A}\,$ with the corresponding Riemann invariants
$R=u-\sqrt{A}\ln \rho$ and $S=u+\sqrt{A}\ln \rho$. Their derivatives, $r$ and $s$, satisfy the
 pair of equations, corresponding to (\ref{eq:gooda}),(\ref{eq:goodb})  with $\theta=(\gamma-1)/2=0$,

\begin{subequations}\label{eq:g1}
\begin{align}\label{eq:g1a}
\disp\bigg(\frac{r}{\srho}\bigg)^{\dac} +\frac12
\disp\frac{r^2}{\srho}&=k\srho, \\
\label{eq:g1b}
 \disp{\bigg({\frac{s}{\srho}}\bigg)}^{\prime}+\frac12
\disp\frac{s^2}{\srho}&=k\srho.
\end{align}
\end{subequations}

\section{Global smooth solutions for sub critical initial data} \label{CT}

\setcounter{figure}{0}
\setcounter{equation}{0}

For the pressure-less Euler-Poisson equations  (\ref{PEP}), the
evolution of $u_x$ and $\rho$ could be traced backwards along the
same particle path to their initial data at $t=0$. The scenario
becomes more complicated with the additional pressure term, due to
the coupling of $r$ and $s$ along \emph{different} particle paths
which are traced back to different neighborhoods of the initial line
$t=0$. This is the main obstacle in finding the sharp critical
threshold of the full Euler-Poisson system (\ref{FEP}). To this end,
we will seek invariant regions for the coupled system, governing the
Riemann invariants. We begin this section with the following lemma.

\begin{lemma} \label{Lemma1}
Given that the total charge $E_0:=\int_{-\infty}^{\infty}\rho_0(x)dx<\infty $\,,
then $\rho(x,t)$ and $u(x,t)$ remain uniformly bounded for all $t>0$.
\end{lemma}
\begin{proof} Under the given condition, we can set (e.g., \cite[p. 116]{ELT01})
$$\varphi_x(x,t)=\frac12\bigg(\int_{-\infty}^x
\rho(\xi,t)d\xi-\int_x^{\infty}\rho(\xi,t)d\xi\bigg)\,,$$
which satisfies $-E_0\leq \varphi_x(x,t)\leq E_0$\,, for all $t\geq 0$
and $x\in \mathbb{R}$\,.\\
Recall the transport equations (\ref{RIa}),(\ref{RIb}) which govern the Riemann invariants
along different characteristics
$R^{\backprime}+k\varphi_x= S^{\prime}+k\varphi_x=0$.
Since $\varphi_x$ is bounded, these transport equations tell us that $R$ and $S$ remain uniformly bounded with at most a linear growth in time. Indeed, for all $M \gg 1$ we have
\begin{equation}\label{eq:sRS}
\sup_{|x|\leq M} \left\{|R(x,t)|, |S(x,t)|\right\} \leq C_{0} + kE_0t, \quad
C_{0}:=\sup_{ |x|\leq M+u_\infty t}\left\{|R_0(x)|, |S_0(x)|\right\}.
\end{equation}
Take the sum and difference of $S$ and
$R$ to find that $u(x,t)$ and $\rho(x,t)$ in (\ref{eq:RS}) remain bounded,

\begin{align}\label{eq:surho}
u_\infty:=\sup_{|x|\leq M} |u(x,t)|   \leq C_{0} + kE_0t, \qquad
\sup_{|x|\leq M} \rho(x,t)  \leq Const. \left\{
\begin{array}{ll} {\displaystyle \left(C_{0} + kE_0t\right)^{\frac{2}{\gamma-1}}}, & \gamma >1, \\
                  {\displaystyle \exp\left(kE_0t\right)},   &  \gamma=1.
\end{array}
\right.
\end{align}
\end{proof}

\noindent
We note in passing that the time growth asserted in (\ref{eq:surho}) is probably  not sharp; the estimate can be improved after taking into account the \emph{uniform bounds} of $R_x/\sqrt{\rho}$ and $S_x/\sqrt{\rho}$ discussed in theorems \ref{MT1} and \ref{MT} below.

\begin{nota}
According to Lemma \ref{Lemma1}, the only way that the full
Euler-Poisson system (\ref{FEP}) breaks down at a finite time is
through the formation of shock discontinuities where $|u_x|$ and/or
$|\rho_x|$ blow up  $\uparrow\infty$\,, but neither will concentrate
at any critical point. This is in contrast to the breakdown of the
``pressure-less" Euler-Poisson equations (\ref{PEP}), where
$-u_x(x,t)=\rho(x,t)\uparrow\infty$ simultaneously at the critical
time.
\end{nota}

\subsection{Critical Threshold for isothermal case: $\gamma=1$}\label{subsec:31}
We begin with the isothermal case, $\gamma=1$, which plays an important
role in various applications. Compared with the general case
(\ref{pair:good}), the isothermal case  becomes simpler due to
the fact that $\theta=0$  decouples the dependence on $r$ and $s$
through the mixed term $\theta rs$, which disappears from left hand
side of (\ref{eq:g1}). Here we  prove the following  sharp
characterization of the critical threshold  phenomena.

\begin{theorem}\label{MT1}
Consider the isothermal Euler-Poisson system (\ref{FEP}) with pressure forcing
$p(\rho)=A\rho$\,, and subject to  initial data $(u_0,\rho_0>0)$ with
finite total charge, $E_0=\int_{-\infty}^{\infty}\rho_0(x)dx<\infty$.
The system admits a global smooth, $C^1$-solution if and only if
\begin{equation}\label{eq:CTgeq1}
u_{0x}(x)\geq -\sqrt{2k\rho_0(x)}+
\disp\sqrt{A}\frac{|\rho_{0x}(x)|}{\rho_0(x)}\,, \quad \forall x \in
\mathbb{R}\,.
\end{equation}
\end{theorem}
\begin{nota}
Expressed in terms of the Riemann invariants specified in
 \S\ref{subsec:iso}, $u_x\pm \sqrt{A}\rho_x/\rho$, theorem \ref{MT1}
states that the isothermal Euler-Poisson equations admit  global
smooth solutions for sub-critical initial conditions,
\begin{equation}\label{eq:rsrho}
s_0\geq -\sqrt{2k\rho_0} \quad {\rm and}  \quad r_0\geq
-\sqrt{2k\rho_0}.
\end{equation}
\end{nota}
\begin{proof} We define $X:=\disp\frac{r}{\srho}$ and $Y:=\disp\frac{s}{\srho}$.
Equations (\ref{eq:g1a}),(\ref{eq:g1b}) then read

\begin{subequations}
\begin{align}
\label{X} \disp X ^{\dac} =\frac{\srho}{2}(2k-X^2), \\
\label{Y} \disp {Y}^{\prime}=\frac{\srho}{2}(2k-Y^2).
\end{align}
\end{subequations}
It follows that
$$ X^{\dac}
\left\{ \begin{array}{ll} >0\,, & X\in(-\sqrt{2k},\sqrt{2k}\,)\,,\\
=0\,, & |X|=\sqrt{2k} \,, \\ <0\,, & |X|>\sqrt{2k},
\end{array} \right. $$
and similarly,
$$ Y^{\prime}
\left\{ \begin{array}{ll} >0\,, & Y \in(-\sqrt{2k},\sqrt{2k}\,)\,,\\
=0\,, & |Y| = \sqrt{2k}\,, \\ <0\,, & |Y| > \sqrt{2k}\,.
\end{array} \right. $$
Thus, starting with (\ref{eq:rsrho}), $X_0,Y_0\geq -\sqrt{2k}$,
we find that $X$ and $Y$ remain bounded within the invariant region $[-\sqrt{2k}, \sqrt{2k}]$, or otherwise, they are decreasing outside this interval. We conclude that
\[
X(\cdot,t),Y(\cdot,t) \leq \max\left\{\sqrt{2k}, X_0(\cdot),Y_0(\cdot)\right\}.
\]
Lemma \ref{Lemma1} tells us that $\rho$ is bounded. The boundedness
of $X$\,, $Y$ and $\rho$ imply that $r=X\srho$ and $s=Y\srho$ remain
bounded for all $t< \infty$, and hence the Euler-Poisson system
(\ref{FEP}) admits a global smooth $C^1$-solution.

Conversely, suppose that there exists $X_0=X(x_0)<-\sqrt{2k}$. We
will show that this value will evolve along
$\Gamma_{\lambda}(x_0,0)$ such that $X(\cdot,t)$ will tend to
$-\infty$  at a finite time. To this end, assume that $Y$ is well
behaved, i.e., $Y_0(\cdot)\geq -\sqrt{2k}$ so  that $Y(\cdot,t)\leq
Y_1:=\max\left\{Y_0(\cdot),\sqrt{2k}\right\}$ for all $t$'s
(otherwise, the finite time blow up of $Y$ can be argued along the
same lines). It follows that $s=Y\srho\leq Y_1\srho$ and inserting
this into $\rho^{\backprime}=-\rho s$, we find
$\rho^{\backprime}\geq-Y_1\rho^{3/2}$. This yields the lower-bound
\[
\rho \geq \bigg(\disp \frac{2}{Y_1t+2/\sqrt{\rho_0}}\bigg)^2,
\]
and together with (\ref{X}), we conclude that $X(\cdot,t)$ satisfies
the following Ricatti equation along the $\Gamma_\lambda$-path,
\begin{equation}\label{Xd}
X^{\backprime}\leq -\disp\frac{X_1}{Y_1t+2/\sqrt{\rho_0}}X^2, \qquad
X_1:=(X_0^2-2k)/X_0^2 >0.
\end{equation}
Integration of (\ref{Xd}) yields
\begin{equation}
X(\cdot,t)\leq \frac{Y_1}{X_1\ln\left(1+\sqrt{\rho_0}Y_1t/2\right)+Y_1X_0}
\end{equation}
Thus, starting with $X_0<-\sqrt{2k}<0$ it follows that there exists
a finite critical time $t_c>0$ such that $X(t\uparrow t_c)$ tends to
$-\infty$.
\end{proof}

The critical threshold condition (\ref{eq:CTgeq1}) reflects the competition between the Poisson
forcing and the pressure. It yields global smooth solutions for a ``large" set of initial configurations allowing negative velocity gradients.
In the particular case that there is no pressure, $A=0$\,, (\ref{eq:CTgeq1}) is reduced to the critical threshold condition of the ``pressure-less" Euler-Poisson equations $u_{0x}> -\sqrt{2k\rho_0(x)}$ of \cite{ELT01}.


\subsection{Critical threshold for $\gamma >1$}\label{sec32}
The equations for the Riemann invariants (\ref{eq:gooda}),
(\ref{eq:goodb}) are coupled through the mixed term, $\theta rs/2$.
We note in passing that it is possible to get rid of this mixed term
when integrating  (\ref{sys}), (\ref{sysa}) with  the integration
factors  $\rho^{(\gamma-3)/4}$\,, and  $
r\rho^{(\gamma-7)/4}(3-\gamma)/4$ in the first and second equations
in each pair, yielding

\begin{align*}
\disp\bigg({r}{\rho^{\frac{\theta-1}{2}}}\bigg)^{\dac} +
\frac{1+\theta}{2}\disp{r^2}{\rho^{\frac{\theta-1}{2}}}
& =k\rho^{\frac{1+\theta}{2}}, \\
\disp{\bigg({{s}{\rho^{\frac{\theta-1}{2}}}}\bigg)}^{\prime}
+\frac{1+\theta}{2} {s^2}{\rho^{\frac{\theta-1}{2}}} &
=k\rho^{\frac{1+\theta}{2}}.
\end{align*}

Nevertheless, it will prove useful to use the same integration
factors, $1/\srho$ and $ r/2\rho\srho$ which led to
(\ref{pair:good}). The main task is to identify the invariant region
associated with (\ref{pair:good}), corresponding to the isothermal
invariant region $[-\sqrt{2k},\sqrt{2k}]$ discussed in theorem
\ref{MT1}.

\begin{theorem}\label{MT}
Consider the Euler-Poisson system (\ref{FEP}) with $\gamma$ law
pressure $p(\rho)=A\rho^{\gamma}$\,, $\gamma>1$, subject to initial
data $(u_0,\rho_0>0)$ with finite total charge,
$E_0=\int_\infty^\infty \rho_0(x)dx < \infty$. Then, there exists a
constant $K_0>0$ depending on $k, \gamma$ and the initial conditions
(specified in (\ref{eq:ctrsb}) below), such that the Euler-Poisson
equations (\ref{FEP}) admit a global smooth, $C^1$-solution if,
\begin{equation}\label{ctg}
u_{0x}(x) \geq -K_0\sqrt{\rho_0(x)} +
\sqrt{A\gamma}\frac{|\rho_{0x}(x)|}{\rho_0(x)^{\frac{3-\gamma}{2}}}.
\end{equation}
\end{theorem}

\noindent
Before we turn to the proof of this theorem, several remarks are in order.

\begin{nota}
Expressed in terms of the Riemann invariants, $r=u_x-
\sqrt{A\gamma}\rho_{0x}/\rho_0^{(3-\gamma)/2}$ and $s=u_x+
\sqrt{A\gamma}\rho_{0x}/\rho_0^{(3-\gamma)/2}$, the critical
threshold (\ref{ctg}) reads
\begin{subequations}\label{eq:ctrs}
\begin{equation}\label{eq:ctrsa}
\frac{r_0(x)}{\sqrt{\rho_0(x)}}, \frac{s_0(x)}{\sqrt{\rho_0(x)}} \geq -K_0.
\end{equation}
The constant $K_0$ is given by
\begin{equation}\label{eq:ctrsb}
K_0=\frac{-\theta M_0+\sqrt{\theta^2 M_0^2 +
8k(1+\theta)}}{2(1+\theta)}, \qquad M_0=\max_x \bigg\{\sqrt{2k},
\frac{r_0(x)}{\sqrt{\rho_0(x)}}, \frac{s_0(x)}{\sqrt{\rho_0(x)}}
\bigg\}
\end{equation}
\end{subequations}
\end{nota}

We mention two simplifications which are summarized in the following
two corollaries. We first observe that if the initial configurations
satisfy the \emph{upper-bound} $r_0(x),s_0(x) \leq \sqrt{2k
\rho_0(x)}$ then (\ref{eq:ctrsb}) yields ${\displaystyle
M_0=\sqrt{2k}}$, hence ${\displaystyle
K_0=\frac{\sqrt{2k}}{1+\theta}}$, and  theorem \ref{MT} implies the
following.

\begin{cor}\label{cor:mta}
Consider the Euler-Poisson system (\ref{FEP}) with $\gamma$ law
pressure $p(\rho)=A\rho^{\gamma}$\,, $\gamma >1$, subject to initial
data $(u_0,\rho_0>0)$ with finite total charge,
$E_0=\int_\infty^\infty \rho_0(x)dx < \infty$. Then, the
Euler-Poisson equations (\ref{FEP}) admit a global smooth,
$C^1$-solution if for all $x\in \R$,
\begin{equation}\label{cta}
|u_{0x}(x)| \leq \sqrt{2k \rho_0(x)} - \sqrt{A\gamma}\frac{|\rho_{0x}(x)|}{\rho_0(x)^{\frac{3-\gamma}{2}}}.
\end{equation}
\end{cor}

\noindent The next result follows from the trivial inequality
${\displaystyle -K_0 \leq \frac{\theta M_0-\big(\theta M_0 +
\sqrt{8k(1+\theta)}\big)/\sqrt{2}}{2(1+\theta)}}$.

\begin{cor}\label{cor:mtb}
Consider the Euler-Poisson system (\ref{FEP}) with a $\gamma$-law
pressure $p(\rho)=A\rho^{\gamma}$\,, $\gamma >1$, subject to initial
data $(u_0,\rho_0>0)$ with finite total charge,
$E_0=\int_\infty^\infty \rho_0(x)dx < \infty$. Then, the
Euler-Poisson equations (\ref{FEP}) admit a global smooth,
$C^1$-solution, if for all $x\in \R$,

\begin{align}
u_{0x}(x)  \geq & -\sqrt{\frac{2k \rho_0(x)}{\gamma+1}} + \nonumber \\
 & + \Big(1-\frac{1}{\sqrt{2}}\Big)\frac{\gamma-1}{2(\gamma+1)}\max_x \bigg\{\sqrt{{2k \rho_0(x)}},
u_{0x}(x) +\sqrt{A\gamma} \frac{|\rho_{0x}(x)|}{\rho_0(x)^{\frac{3-\gamma}{2}}}\bigg\}
+\sqrt{A\gamma}\frac{|\rho_{0x}(x)|}{\rho_0(x)^{\frac{3-\gamma}{2}}}. \label{ctb}
\end{align}
\end{cor}

\begin{nota}
We observe that as in the isothermal case, the  critical threshold
in its various versions (\ref{ctg}),(\ref{eq:ctrs}), (\ref{cta}) and
(\ref{ctb}),  allow a ``large set" of initial configurations with
negative velocity gradient, due to the competition between the
stabilizing Poisson forcing $k\rho\phi(\rho)_x$ and the
destabilizing pressure $A(\rho^\gamma)_x$. In the extreme case that
Poisson forcing is missing $k=0$\,, the breakdown of the system is
generic
unless $u_{0x}$ is   
\emph{positive enough} (so that $r_0,s_0>0$). In the other extreme
of a ``pressure-less" Euler-Poisson, $A\!=\!0, \gamma\!=\!1$\,, the
critical thresholds (\ref{ctg}), (\ref{cta}) are reduced to
$u_{0x}(x)>-\sqrt{2k\rho_0(x)}$\,, which coincides with the
``pressure-less" critical threshold  (\ref{epct}) found in
\cite{ELT01}.
\end{nota}

\begin{proof} Expressed in terms of $\disp X:=\frac{r}{\srho}$ and $Y:=\disp \frac{s}{\srho}$, equations
(\ref{pair:good}) read
\begin{subequations}\label{eq:XY}
\begin{align}\label{eq:X}
X^{\dac} = \srho\bigg(k-\frac{1+\theta}{2}X^2 + \frac{\theta}{2} XY\bigg),\\
\label{eq:Y}
Y^{\prime}= \srho\bigg(k-\frac{1+\theta}{2}Y^2 +
\frac{\theta}{2} XY\bigg).
\end{align}
\end{subequations}
We seek an invariant region of the form $[-K_0,M_0]$, with $K_0,M_0>0$ yet to be determined.
To this end we construct a  ``buffer zone" in which  positive values of $X,Y$ must decrease and hence remain upper-bounded. We begin by noticing that if $X,Y\leq  M$ then\footnote{We let $Z_+=\max\{X, 0\}$ and
 $Z_-=\min\{Z,0\}$ denote the positive and negative part of $Z$.}
$X_+Y \leq  M^2$, and recalling that $\theta\geq0$, (\ref{eq:XY}) then yields
\begin{subequations}\label{eq:XYM}
\begin{align}
X^{\dac} \leq  \srho\bigg(k-\frac{1+\theta}{2}X^2 + \frac{\theta}{2} M^2\bigg), \quad X>0,\\
Y^{\prime} \leq  \srho\bigg(k-\frac{1+\theta}{2}Y^2 +
\frac{\theta}{2} M^2\bigg), \quad Y>0.
\end{align}
\end{subequations}
This in turn implies that
\[
X \ {\rm or } \ Y \ {\rm \ are \ decreasing \ whenever}  \ X \in  {\mathcal I}_M \ {\rm or} \  Y  \in  {\mathcal I}_M.
\]
Here ${\mathcal I}_M$ is the interval ${\mathcal I}_M:=\big(C_+(M), M\big)$ where 
$C_+(M):= \sqrt{(2k+\theta M^2)/(1+\theta)}$ is dictated by the largest root of the quadratics on the right of (\ref{eq:XYM}). Observe that for ${\mathcal I}_M$ to be nonempty requires $M>\sqrt{2k}$. We therefore set, 
$M_\eps:=\max_x\bigg\{\sqrt{{2k}}+\eps, X_0(x), Y_0(x) \bigg\}$. 
We claim that $X,Y \leq M_\eps$: indeed, either $X, Y \leq C_+(M_\eps) < M_\eps$ or, if  $X,Y > C_+(M_\eps)$, then they must decrease being ``traped" inside ${\mathcal I}_{M_\eps}$ and hence $X,Y \leq M_\eps$.  Letting $\eps \downarrow 0$ we end up with the upper-bound
\begin{equation}\label{eq:upp}
X(\cdot,t),Y(\cdot,t) \leq M_0, \qquad M_0:=
\max_x\bigg\{\sqrt{{2k}}, X_0(x), Y_0(x) \bigg\}.
\end{equation}
In a similar manner, we study the lower bound of the invariant region. By (\ref{eq:upp})
 and (\ref{eq:XY}) yield
\begin{subequations}\label{sub:stam}
\begin{align}
X^{\dac} \geq  \srho\bigg(k-\frac{1+\theta}{2}X^2 + \frac{\theta}{2} M_0X\bigg), \quad X<0,\\
Y^{\prime} \geq  \srho\bigg(k-\frac{1+\theta}{2}Y^2 +
\frac{\theta}{2} M_0Y\bigg), \quad Y<0,
\end{align}
\end{subequations}
which in turn, imply that
\begin{subequations}
\begin{equation}\label{eq:XK}
X \ {\rm and} \ Y \ {\rm \ are \ increasing \ if \ }  \ 0 \geq X,Y > -K_0,
\end{equation}
where $K_0$ is the smallest root of the quadratics on the right of (\ref{sub:stam}),
\begin{equation}\label{eq:XKb}
 K_0:= \frac{-\theta M_0+\sqrt{\theta^2M_0^2+8k(1+\theta)}}{2(1+\theta)}.
\end{equation}
\end{subequations}
 The critical threshold condition (\ref{ctg}) tells us that at $t=0$,
$X_0,Y_0 \geq -K_0$ and (\ref{eq:XK}) implies that $X(\cdot,t)$ and
$Y(\cdot,t)$ remain  above the same lower-bound, (\ref{ctg}). As
before, the bounds of $X, Y$ and $\rho$ imply that $r=X\srho$ and
$s=Y\srho$ remain bounded, and hence the Euler-Poisson system
(\ref{FEP}) a global smooth, $C^1$-solution.
\end{proof}

\section{Finite time breakdown for super-critical initial data} \label{FTB}
\setcounter{figure}{0}
\setcounter{equation}{0}
Consider the Euler-Poisson system (\ref{FEP}) with a $\gamma$-law
pressure, $\gamma \geq 1$, and subject to
initial data such that  $r_0(x), s_0(x) \leq \sqrt{2k}$. Then, according to corollary \ref{cor:mta}, the following critical threshold is sufficient for the existence of global smooth solutions,
\[
u_{0x}(x) \geq -\sqrt{2k \rho_0(x)} + \sqrt{A\gamma}\frac{|\rho_{0x}(x)|}{\rho_0(x)^{\frac{3-\gamma}{2}}}.
\]
In this section we show that this critical threshold is also \emph{necessary} for global regularity.
\begin{theorem}\label{BD}
Consider the Euler-Poisson system (\ref{FEP}) with a $\gamma$-law
pressure $p(\rho)=A\rho^{\gamma}$\,, $\gamma \geq 1$, subject to
initial data $(u_0,\rho_0>0)$\,. The system loses the
$C^1$-smoothness if there exists an $x\in \mathbb{R}$ such that
\begin{equation}\label{bdc}
u_{0x}(x) < -\sqrt{2k\rho_0(x)} +
\sqrt{A\gamma}\frac{|\rho_{0x}(x)|}{\rho_0(x)^{\frac{3-\gamma}{2}}}.
\end{equation}
\end{theorem}

\begin{nota}
Expressed in terms of the Riemann invariants, $r=u_x-
\sqrt{A\gamma}\rho_{0x}/\rho_0^{(3-\gamma)/2}$ and $s=u_x+
\sqrt{A\gamma}\rho_{0x}/\rho_0^{(3-\gamma)/2}$, the condition
(\ref{bdc}) reads
\begin{equation}\label{bdcri}
\exists\, x\in \mathbb{R}\quad \rm{s.t.} \quad
r_0(x)<-\sqrt{2k\rho_0(x)}\,, \quad \rm{or} \quad
s_0(x)<-\sqrt{2k\rho_0(x)}\,.
\end{equation}
The lack of smoothness in this case was shown in theorem \ref{MT1} for $\gamma=1$ and is extended for $\gamma >1$ below.
\end{nota}

\begin{proof}
Recall equations (\ref{eq:XY}) for $\disp X:=\frac{r}{\srho}$ and
$Y:=\disp \frac{s}{\srho}$
\begin{subequations}
\begin{align}
X^{\dac} = \srho\bigg(k-\frac{1+\theta}{2}X^2 + \frac{\theta}{2} XY\bigg),\\
Y^{\prime}= \srho\bigg(k-\frac{1+\theta}{2}Y^2 + \frac{\theta}{2}
XY\bigg).
\end{align}
\end{subequations}
In the proof of theorem \ref{MT}, we have shown that $X$ and $Y$
have an upper bound
\begin{equation}
X(\cdot,t),Y(\cdot,t) \leq M_0, \qquad M_0:=
\max_x\bigg\{\sqrt{{2k}}, X_0(x), Y_0(x) \bigg\}.
\end{equation}
Suppose that there exists $X_0=X(x_0)<-\sqrt{2k}$. We will show that
this value will evolve along $\Gamma_{\lambda}(x_0,0)$ such that
$X(\cdot,t)$ will tend to $-\infty$  at a finite time. To this end,
assume that $Y$ is well behaved, i.e., $Y_0(\cdot)\geq -\sqrt{2k}$
so  that $Y(\cdot,t)\leq M_0$ for all $t$'s (otherwise, the finite
time blow up of $Y$ can be argued along the same lines). It follows
that along $\Gamma_{\lambda}(x_0,0)$
\begin{equation}
X^{\dac} = \srho\bigg(k-\frac{1+\theta}{2}X^2 + \frac{\theta}{2}
XY\bigg) < \srho\bigg(k-\frac{1}{2}X^2\bigg)\,.
\end{equation}
Following exactly what we have done in the proof of theorem
\ref{MT1}, we obtain the inequality
\begin{equation}
X(\cdot,t)\leq
\frac{M_0}{X_1\ln\left(1+\sqrt{\rho_0}M_0t/2\right)+M_0X_0}\,,
\end{equation}
where $X_1:=(X_0^2-2k)/X_0^2 >0$\,. Thus, starting with
$X_0<-\sqrt{2k}<0$ it follows that there exists a finite critical
time $t_c>0$ such that $X(t\uparrow t_c)$ tends to $-\infty$.
\end{proof}

\noindent
We conclude with an example for a finite time breakdown.

\noindent {\bf Example}: Suppose at $t=0$\,, $u_0(x)=0$\, and
$$\rho_0(x)=\left\{ \begin{array}{ll} 1\,,  & x<0\,, \\
 1-\disp\frac{x}{2\epsilon}\,, & 0\leq x \leq \epsilon\,,  \\ \frac12\,,   &
 x>\epsilon\,.
\end{array}\right.$$
 Thus $$s_0(x)=\left\{
\begin{array}{ll} -\sqrt{A\gamma}\Big(1-\disp\frac{x}{2\epsilon}\Big)/2\epsilon\,,
& 0<x<\epsilon\,, \\ 0\,, & \rm{elsewhere}\,.
\end{array} \right. $$
If we choose $\epsilon$ small enough, then
$s_0(x)<-\sqrt{2k\rho_0(x)}$ for $0<x<\epsilon$\,. According to
theorem \ref{BD}, the system (\ref{FEP}) will break down in a finite
time. This example shows that even if the fluid is near rest at
$t=0$, the pressure itself could still lead to collision.


\end{document}